\documentclass[12pt]{article}
\usepackage[english]{babel}
\usepackage[centertags]{amsmath}
\usepackage{amsfonts}
\usepackage{amssymb}
\usepackage{amsthm}
\usepackage{newlfont}
\usepackage[bookmarksnumbered, colorlinks, plainpages]{hyperref}

\usepackage{natbib}
\begin{document}

\title{An intrinsic characterization of semi-normal operators}
\author{L. Z. Gevorgyan\\ State Engineering University of
Armenia,\\ Department of Mathematics, E-mail: levgev@hotmail.com}
\date{}
\maketitle

AMS MSC {47A12,\{Primary\}, 51N15, \{Secondary\}}

\begin{abstract}
Two necessary and sufficient conditions for an operator to be semi-normal are revealed. For a Volterra integration operator the set where the operator and its adjoint are metrically equal is described.
\end{abstract}

Let  $A$ be a linear bounded operator, acting in a Hilbert
space $\left ( {\mathcal{H},\left\langle { \bullet , \bullet }
\right\rangle }\right)$  and $W\left(A\right)$ denote  the numerical range of $A.$ If $C\left(A \right)=A^*A  - AA^*$ is semi-definite, the operator $A$ is said \cite{MR0217618} to be semi-normal, particularly, if $C\left( A \right) \ge \textbf{0},$
then $A$ is hyponormal. The well-known and important class of normal operators is characterized by the equality $ AA^*  = A^*A.$ It is easy to see that the last condition is equivalent to the equality $ \left\| {Ax} \right\| = \left\| {A^* x} \right\|$ for any $x\in \mathcal{H},$ meaning that any normal operator is metrically equal to its adjoint on all $\mathcal{H}.$ For hyponormal operator in \cite{MR0187097} is proved that conditions \begin{equation} \left\| {Ax} \right\| = \left\| {A^* x} \right\|\,\, \text{and}\,\,\,\   A^* Ax = AA^* x \label{Gevorgyaneq1}\end{equation} are equivalent. Note that the set of points, satisfying the second condition is the null space of the self-commutator- $N\left( C\left(A \right)\right).$
As the both conditions are symmetric, Stampfli's result remains valid for semi-normal operators. Using this property, Stampfli has shown that any extreme point of the numerical range of a hyponormal operator $A$ is a reducing eigenvalue.

Denote \[E\left( A \right) = {\text{ }}\left\{ {x:\left\| {Ax} \right\| = \left\| {A^* x} \right\|{\text{ }}} \right\}\]
and \[M_\lambda  \left( A \right) = \left\{ {x:\left\langle {Ax,x} \right\rangle  = \lambda \left\| x \right\|^2 } \right\}\]
Evidently conditions $\lambda  \in W\left( A \right)$ and $M_\lambda  \left( A \right) \ne \left\{ \theta  \right\}$ are equivalent.

\emph{\textbf{Proposition 1.}} For any operator $A$ one has \[
\left\| {Ax} \right\|^2  - \left\| {A^* x} \right\|^2  = \left\langle {C\left( A \right)x,x} \right\rangle ,
\]
particularly, \[E\left( A \right) = M_0 \left( C\left(A \right) \right).\]

\emph{\textbf{Proof.}} As $\left\| {Ax} \right\|^2  = \left\langle {Ax,Ax} \right\rangle  = \left\langle {A^* Ax,x} \right\rangle$ and
$\left\| {A^* x} \right\|^2  = \left\langle {AA^* x,x} \right\rangle,$  the conditions
$\left\| {Ax} \right\| = \left\| {A^* x} \right\|$ and $\left\langle {\left( C\left(A \right) \right)x,x} \right\rangle  = 0$ are equivalent.

\emph{\textbf{Proposition 2.}} The operator $A$ is semi-normal if and only if $0$ is an extreme point of the closure of
$W\left( C\left(A \right) \right).$

\emph{\textbf{Proof.}} As the numerical range of any self-adjoint operator is a convex subset of $\mathbb{R}$ we have $\overline {W\left( {C\left( A \right)} \right)}  = \left[ {a;b} \right].$ If  $0$ is an extreme point of $\overline {W\left( {C\left( A \right)} \right)}$ then $ab=0,$ hence $A$ is semi-normal. Let now $A$ be semi-normal, i.e $ab\ge 0.$ We show that the strict inequality $ab > 0$ is not possible. According to a result of Radjavi (\cite{MR0203482}, Corollary 1) if $B$ is  a selfadjoint operator such that $B \geqslant \alpha I\,\,\left( {B \leqslant -\alpha I} \right)$ for some positive number $\alpha,$ then $B$ is not a self-commutator.Thus $ab=0$ and $0 $ is an extreme point of $\overline {W\left( {C\left( A \right)} \right)}.$

\emph{\textbf{Proportion 3.}} The equivalence \eqref{Gevorgyaneq1} is true if and only if the operator $A$ is semi-normal.

\emph{\textbf{Proof.}} Only the necessity of this condition should be proved. Let \eqref{Gevorgyaneq1} be true. If $E\left( A \right) = \left\{ \theta  \right\},$ then $0 \notin W\left( C\left(A \right) \right),$ hence it lies entirely in the positive or negative semi-axis.
Let now $x$ and $y$ be two elements from $E\left( A \right).$ Then from
$\left\| {Ax} \right\| = \left\| {A^* x} \right\|,\left\| {Ay} \right\| = \left\| {A^* y} \right\|$ follows
$AA^* x = A^* Ax,AA^* y = A^* Ay$ and $AA^* \left( {x + y} \right) = A^* A\left( {x + y} \right),$ implying $\left\| {A\left( {x + y} \right)} \right\| = \left\| {A^* \left( {x + y} \right)} \right\|.$
According to \cite{MR0259634} the linearity of $M_\lambda  \left( A \right)$ is equivalent to the condition that $\lambda$ is an extreme point of $W\left(A\right).$ Thus $0$ is an extreme point of $W\left(C\left(A \right) \right),$ completing the proof.

\emph{\textbf{Remark 1.}} The principal reason in the proof above was the linearity of $M_0\left(C\left(A \right) \right).$ If the last condition is satisfied, then  $A$ is semi-normal and by Stampfli's result $E\left(  A \right) = N\left( C\left(A \right) \right).$

\emph{\textbf{Remark 2.}} In (\citet{MR2320391}, Proposition 2, Corollary 1) is proved that $
N\left( A \right) = M_0 \left( A \right)$ if and only if \[
A = \left( {\begin{array}{*{20}c}
   B & 0  \\
   0 & 0  \\

 \end{array} } \right),\]  $0 \notin W\left( B \right),$ where any direct summand may be absent.

The situation is more interesting for non semi-normal operators. The example below  exhibits the matter for a non semi-normal quasinilpotent compact operator.

\emph{\textbf{Example.}} Consider the Volterra integration operator $V$
\[\left( {Vf} \right)\left( x \right) = \int\limits_0^x {f\left( t \right)} dt, f\in L^2 \left( {0\,;\,1} \right).\]
We have $V1 = x,V^* 1 = 1 - x,$ implying $\left\| {V1} \right\| = \left\| {V^* 1} \right\|.$ Let now $f \bot 1.$ As $\int\limits_0^x {f\left( t \right)} dt + \int\limits_x^1 {f\left( t \right)} d = \int\limits_0^1 {f\left( t \right)} dt,$ we have $Vf =  - V^* f$ and $\left\| {Vf} \right\| = \left\| {V^* f} \right\|,$ therefore $\left\{ {1,L^2 \left( {0;1} \right)\ominus 1} \right\} \subset E\left( A \right).$

The self-commutator of $V$ is \[
\left( {C\left( V \right)f} \right)\left( x \right) = \int\limits_0^1 {f\left( t \right)} dt - x\int\limits_0^1 {f\left( t \right)} dt - \int\limits_0^1 {tf\left( t \right)} dt\]
or
\[\left( {C\left( V \right)f} \right)\left( x \right) = \left( {\frac{1}
{2} - x} \right)\int\limits_0^1 {f\left( t \right)} dt + \int\limits_0^1 {\left( {\frac{1}
{2} - t} \right)f\left( t \right)} dt.\]
Denoting $e_1  = 1,\,e_2  = \sqrt 3 \left( {1-2x} \right)$ (they are two first $L^2 \left( {0;1} \right)$-orthonormal polynomials)
we get \[
C\left( V \right)f = \frac{1}
{{2\sqrt 3 }}\left( {\left\langle {f,e_1 } \right\rangle e_2  + \left\langle {f,e_2 } \right\rangle e_1 } \right).\] Now we introduce two orthonormal elements \[
\begin{gathered}
  u_1  = \frac{1}
{{\sqrt 2 }}\left( {e_1  + e_2 } \right) = \sqrt {2 + \sqrt 3 }  - \sqrt 6 x, \hfill \\
  u_2  = \frac{1}
{{\sqrt 2 }}\left( {e_1  - e_2 } \right) = \sqrt 6 x - \sqrt {2 - \sqrt 3 } , \hfill \\
\end{gathered}
\]
and arrive at the canonical form of the self-commutator of $V$ \begin{equation}
C\left( V \right)f = \frac{1}
{{2\sqrt 3 }}\left( {\left\langle {f,u_1 } \right\rangle u_1  - \left\langle {f,u_2 } \right\rangle u_2 } \right).\label{canon}\end{equation}
Note that the product $u_1 u_2$ defines  the third orthogonal polynomial $6x^2  - 6x + 1.$

From \eqref{canon} follows that the spectrum of $C\left( V \right)$ is the set $\left\{ { - \frac{{\sqrt 6 }}
{2},0,\frac{{\sqrt 6 }}
{2}} \right\},$ hence $
W\left( {C\left( V \right)} \right) = \left[ { - \frac{{\sqrt 6 }}
{2};\frac{{\sqrt 6 }}
{2}} \right].$
The null-space of $C\left( V \right)$ consists of functions orthogonal to the first-order polynomials-
$L^2 \left( {0;1} \right) \ominus  \bigvee \left\{ {1,x} \right\},$ where $\bigvee$ denotes the linear span of the set.
As \[
\left\langle {C\left( V \right)f,f} \right\rangle  = \frac{1}
{{2\sqrt 3 }}\left( {\left| {\left\langle {f,u_1 } \right\rangle } \right|^2  - \left| {\left\langle {f,u_2 } \right\rangle } \right|^2 } \right),
\]
we get
$E\left( V \right) = \left\{ {f:\left| {\left\langle {f,u_1 } \right\rangle } \right| = \left| {\left\langle {f,u_2 } \right\rangle } \right|} \right\},$ i.e. \[
E\left( V \right) = \bigcup\limits_{\varphi  \in \left[ {0;2\pi } \right)} {L_\varphi  }
\]
where $L_\varphi$ is the orthocomplement to the subspace, generated by the element $u_1  - e^{i\varphi } u_2.$

\bibliography{mybib3}

\end{document}